\documentstyle[12pt,a4]{amsart}

\newcommand{\eqn}[1]{(\ref{#1})}

\newcommand{\eeq}{\end{equation}}
\newcommand{\beql}[1]{\begin{equation}\label{#1}}

\newcommand{\ZZ}{{\mathbb Z}}
\newcommand{\RR}{{\mathbb R}}

\newcommand{\HH}{{\mathbb H}}

\newcommand{\QQ}{{\mathbb Q}}
\newcommand{\CC}{{\mathbb C}}

\newcommand{\bv}{{\bf v}}

\newcommand{\bx}{{\bf x}}

\newcommand{\sB}{{\cal B}}

\newcommand{\sD}{{\cal D}}
\newcommand{\sF}{{\cal F}}

\newcommand{\sM}{{\cal M}}

\newcommand{\sO}{{\cal O}}

\newtheorem{theorem}{Theorem}
\newtheorem{lemma}{Lemma}

\begin{document}
\title[\null]{
The Riemann Hypothesis for Certain Integrals of Eisenstein Series
}
\maketitle
\begin{center}
Jeffrey C. Lagarias \\

Masatoshi Suzuki \\
\smallskip

(November 1, 2005)
\end{center}

\begin{abstract}

This paper studies  the  non-holomorphic
Eisenstein series $E(z, s)$ for  the modular
surface $PSL(2, \ZZ)\backslash \HH$, and 
shows that integration with respect to 
certain non-negative  measures $\mu(z)$ 
gives meromorphic functions $F_{\mu}(s)$ 
that have all their zeros on the line $\Re(s) = \frac{1}{2}$.
For the constant
term $a_0(y,s)$ of the Eisenstein series 
the Riemann hypothesis 
holds for  all values $y \ge 1$,
with at most two exceptional real zeros, 
which occur exactly  for those  
$y > 4\pi e^{-\gamma} = 7.0555+$.
The  Riemann hypothesis holds for
all truncation integrals with truncation parameter $T \ge 1$.
At the value $T=1$ this proves the  Riemann hypothesis
for a zeta function $Z_{2, \QQ}(s)$
recently introduced by Lin Weng, associated
to rank $2$ semistable lattices over $\QQ$.
\end{abstract}

%
%
%
%

\setlength{\baselineskip}{1.0\baselineskip}

\section{Introduction}

Basic objects in the theory of automorphic forms are 
Eisenstein series, whose  Fourier coefficients, 
particularly their constant terms, give information
about  $L$-functions. 
We consider the (completed) 
non-holomorphic Eisenstein series $E^{\ast}(z, s)$ for the modular
group $PSL(2, \ZZ)$, which is given for $z= x+iy \in \HH$
with $y > 0$ and $\Re(s) > 1$ by 
\begin{equation}~\label{101}
E^{\ast}(z, s) := \pi^{-s} \Gamma(s) \left( \frac{1}{2} 
\sum_{(m,n) \in \ZZ^{2} \backslash (0,0)} \frac{y^s}{|mz+ n|^{2s}} \right)
 = \pi^{-s} \Gamma(s) E(z, s)
\end{equation}
or, equivalently,
\begin{equation}~\label{102}
E^{\ast}(z, s) = \pi^{-s} \Gamma(s) \zeta(2s)
\left( \frac{1}{2} \sum_{(c,d)=1} 
\frac {y^s}{|cz+d|^{2s}} \right).
\end{equation}
It is well-known that for fixed $z$, $E^{\ast}(z, s)$ 
meromorphically  continues to the $s$-plane, and
satisfies the functional equation
\begin{equation}~\label{103}
E^{\ast}(z, s)= E^{\ast}(z, 1-s),
\end{equation}
and its only singularities are simple poles at $s=0$ and $s=1$
with residues $-\frac{1}{2}$ and $\frac{1}{2}$, respectively. 
In addition $E^{\ast}(z, s)$ behaves like a modular form of weight $0$,
on $ PSL(2, \ZZ) \backslash \HH$, satisfying
\begin{equation}~\label{104}
E^{\ast} \left( \frac{az+b}{cz+d}\,, s \right) = E^{\ast}(z,s)
 ~~~~~~~\mbox{for}~~
\left[ \begin{array}{cc}
a & b \\
c & d
\end{array} \right] \in PSL(2, \ZZ).
\end{equation}
In particular $E^{\ast}(z+1, s) = E^{\ast}(z,s)$ 
so it has a Fourier expansion
\begin{equation}~\label{105}
E^\ast (z, s) = \sum_{n= - \infty}^{\infty} a_n(y, s) e^{2 \pi i n x},
\end{equation}
with
\begin{equation}~\label{106}
a_n(y,s) = \int_0^1 E^\ast (x+iy, s) e^{2\pi i n x} dx.
\end{equation}

The  non-constant Fourier coefficients ($n \ne 0$) are given by 
\begin{equation}~\label{107}
a_n(y, s) = 2|n|^{s - \frac{1}{2}} \sigma_{1-2s}(|n|) \sqrt{y} 
K_{s - \frac{1}{2}}(2 \pi |n|y),
\end{equation}
in which 
\begin{equation}~\label{108}
\sigma_s(n) = \sum_{d~|~n} d^{s} = 
\prod_{p^{e} ||n}\frac{1 - p^{(e+1)s}}{1 - p^s}, 
\end{equation}
and, for positive real $y$,  the $K$-Bessel function is 
\begin{eqnarray}~\label{109}
K_s(y) &= & \int_{0}^{\infty} e^{-\frac{y}{2}(e^u + e^{-u)}} 
\frac{1}{2}(e^{us} + e^{-us})du \nonumber \\
&=& \frac{1}{2}\int_{-\infty}^{\infty} e^{-\frac{y}{2}(x + \frac{1}{x})}
(x^s + x^{-s}) \frac{dx}{x}.
\end{eqnarray}
In particular, $a_n(y, s) = a_{-n}(y, s)$, and for fixed $y$ 
these are entire functions of $s$.
The constant term ($n=0$) is given by the more complicated
expression 
\begin{equation}~\label{110}
a_0(y, s) = \zeta^{\ast}(2s) y^s + \zeta^{\ast}(2-2s) y^{1-s},
\end{equation}
in which $\zeta^{\ast}(s)$ is the {\em completed zeta function}
\begin{equation}~\label{111}
\zeta^{\ast}(s) := \pi^{-\frac{s}{2}} \Gamma(\frac{s}{2})\zeta(s). 
\end{equation}
The constant term $a_0(y,s)$ is a meromorphic function of $s$, and has
simple poles at $s=0$ and $s=1$ with residues $-\frac{1}{2}$
and $\frac{1}{2}$ independent
of $y$; these account for the poles of $E^{\ast}(z,s)$.

For any complex-valued measure $d\mu$ on
on the modular surface
$PSL(2, \ZZ)\backslash \HH$ of finite mass the function
\begin{equation}~\label{112}
F_{\mu}(s) = \int \int_{PSL(2, \ZZ)\backslash \HH} E^\ast(z, s) d\mu(z)
\end{equation}
inherits the functional equation $F_{\mu}(s) = F_{\mu}(1-s)$
from the Eisenstein series. The Fourier coefficients for fixed
$y$ are obtained by such an integral 
\eqn{112} using a one-dimensional
complex measure $d\mu(z) = e^{2\pi i n x} dx$
supported on a closed horocycle at height $y$.
The functional equation for $E(z,s)$ implies  $a_n(y,s) = a_n(y, 1-s)$.

The object of this paper is to study  certain such
integrals, related to the constant term of the Eisenstein series,
and show conditions under which they  
satisfy the Riemann hypothesis: All zeros of $F_{\mu}(s)$
lie on $\Re(s) = \frac{1}{2}.$ The measures we consider are nonnegative
real-valued measures of finite mass.  

The first integral we consider gives a special case of 
zeta functions recently 
introduced by Lin Weng (\cite{we0}, \cite[Sec. B.4]{we1}, \cite{Weng04}) 
whose general definition uses an 
integral representation  motivated
in part by Arakelov geometry. The rank $n$ zeta function $Z_{n, K}(s)$
of a number field $K$ is given in \cite[Sec. II]{Weng04}
as a Mellin-type integral over
a moduli space of rank $n$ semi-stable lattices; this 
can be reduced to an integral of an Eisenstein series associated
to the group $PSL(n, \ZZ)$ over a certain subset of the allowable
lattices. Weng \cite[Main Theorem A ]{Weng04}
shows that this  zeta function meromorphically 
continues to $\CC$, with  singularities
being simple poles at $s=0$
and $s=1$, and satisfies the functional equation
$Z_{n, K}(s)= Z_{n, K} (1-s)$. 
For the case of rank $1$ lattices over $\QQ$, one
obtains $Z_{1,\QQ}(s) = \zeta^{\ast}(s)$, recovering the
usual completed Riemann zeta function.
For the case of rank $2$ lattices over $\QQ$, the
resulting definition can be simplified  to the following integral
of an Eisenstein series:
\begin{equation}~\label{113}
Z_{2, \QQ}(s) := \int_{\sD_{ss}} E^\ast (z, s) d\mu_{\HH}(z) = 
\int \int_{\sD_{ss}} E^\ast (x+iy, s)\frac{dx dy}{y^2},
\end{equation} 
as we explain in the appendix to this paper.
Here the set 
\begin{equation}~\label{113a}
\sD_{ss} := \{ z=x+iy~|~ -\frac{1}{2} \le x \le \frac{1}{2},~ 0<  y \le 1,~
\mbox{and}~  x^2 +y^2 \ge 1\};~ 
\end{equation}
represents the set of two-dimensional semi-stable lattices. 
This integral
can be proved to converge absolutely for
$0 < \Re(s) < 1$. Using it one can deduce that 
$Z_{2, \QQ}(s)$ extends to  a meromorphic function on $\CC$ 
satisfying the 
functional equation $Z_{2, \QQ}(s) = Z_{2, \QQ}(1-s)$,
whose only singularities
are simple poles at $s=0$ and $s=1$, with residue
at $s=1$ (resp. $s=0$) given by  $c_1= \frac{1}{2}(\frac{\pi}{3}-1)$
(resp.  $-c_1$). We show that this function satisfies the Riemann
hypothesis. \\


\begin{theorem}~\label{th11}
The meromorphic function 
\begin{equation}~\label{115}
Z_{2, \QQ}(s)  := \zeta^{\ast}(2s)\frac{1}{s-1}  - 
\zeta^{\ast}(2-2s)\frac{1}{s}
\end{equation}
has all its zeros on the critical line $\Re(s) = \frac{1}{2}$. 
\end{theorem}

A further result in \S4 shows that all these
zeros are simple zeros. These results 
supply evidence that the zeta functions associated
to semi-stable lattices introduced 
by Weng \cite{Weng04} are of interest and  deserve further
investigation.  Note that for $n \ge 2$ 
the rank $n$  zeta function does not have an Euler product.

We next consider functions given by  integration of the
Eisenstein series against hyperbolic measure over the 
{\em truncation region} 
\begin{equation}~\label{116}
\sD^T :=   \{ z=x+iy~|-\frac{1}{2} \le x \le \frac{1}{2} 
\mbox{and}~ y \ge T \}
\end{equation}
in the upper half-plane $\HH$. These regions have finite
hyperbolic measure, 
 and the associated integral is 
\begin{equation}~\label{117}
I(T,s) := \int \int_{\sD^T} E^\ast (x+iy, s)\frac{dx dy}{y^2}. 
\end{equation} 
It can be shown that the  integral $I(T, s)$  converges absolutely 
for all $T >0$, when $0 < \Re(s) < 1$.
The integral \eqn{117} over  the $x$ variable removes all Fourier
terms but the constant term, yielding the alternate formula
\begin{equation}~\label{119}
I(T,s)= \int_{T}^{\infty} a_0(y, s) \frac{dy}{y^2}
= -\zeta^{\ast}(2s) \frac{T^{s-1}}{s-1} + 
\zeta^{\ast}(2-2s) \frac{T^{-s}}{s},
\end{equation}
valid for all $T > 0$. For each $T$ this  function has at most simple poles
at $s=0$ and $s=1$ and satisfies the functional 
equation $I(T, s) = I(T, 1-s)$; its residues at these points
depend on $T$. 
At  $s=1$ the residue is $-\frac{1}{2}(\frac{\pi}{3} - \frac{1}{T})$
and $I(T,s)$ is entire at the point $T= \frac{3}{\pi} = 0.9548+$.
The integration 
region $\sD^T$ is contained in the standard 
fundamental domain 
\begin{equation}~\label{118}
\sF := \{ z|~ |z| > 1,~   -\frac{1}{2} < \Re(z) \le \frac{1}{2} \}
\end{equation}
of the modular surface, when $T \ge 1$. 
We  obtain the following result. \\


\begin{theorem}~\label{th12}
For  each fixed $T \ge 1$, the meromorphic function
\begin{equation}~\label{120}
I(T, s) = -\zeta^{\ast}(2s)\frac{T^{s-1}}{s-1} 
+ \zeta^{\ast}(2-2s) \frac{T^{-s}}{s}
\end{equation}
has all its zeros in the critical line $\Re(s)= \frac{1}{2}.$ 
\end{theorem}

The hypothesis $T \ge 1$ cannot be relaxed.
It can be shown
that the Riemann hypothesis  fails to hold for all values $0 < T < 1$
by an argument principle method as in
Hejhal\cite[p. 89]{Hej90}. 
In \S4 we establish for all $T \ge 1$ that 
all zeros of $I(T,s)$ are simple , using results of \cite{La04}
(see Theorem~\ref{Nth51}).

Theorem~\ref{th12} in the special case $T=1$ yields
Theorem~\ref{th11}, since we have
\begin{equation}~\label{121}
Z_{2, \QQ}(s) = -\int \int_{\sD^{1}} E^{\ast}(z, s) d \mu_{\HH}(z) =- I_1(s).
\end{equation}
 Zagier \cite[Example 1]{Za81c}
observes the Eisenstein series integral identity
$$
\int \int_{\sF} E^\ast (z, s) d \mu_{\HH}(z) \equiv 0, 
$$
for $0 < \Re(s) < 1$. The identity \eqn{121}
follows by combining  this with the
fact that $\sF$ is partitioned into 
the union of $D_{ss}$ and $D^{1}$ 
(up to  a hyperbolic measure zero set).

The second set of  integrals we consider are those 
giving  the constant term
$a_0(y,s)$ for fixed $y$, i.e. \eqn{106} for $n=0$.
Study of these integrals is motivated  by  an observation of Dan Bump
which is stated in \cite[p. 6]{BCKV00}: {\em for each $n \ne 0$ and 
each $y>0$, the Fourier coefficient
$a_n(y,s)$ satisfies the
Riemann hypothesis in the $s$-variable.}
(This fact was noted
earlier by D. Hejhal \cite[p. 85]{Hej90}. Bump's
point in  \cite{BCKV00}
is to make an analogy with metaplectic
Eisenstein series, see the discussion in \S5.)
The observation follows because the finite Dirichlet series
$\sigma_s(n)$ is easily shown from its
Euler product to have all zeros on
the imaginary axis $\Re(s) =0$. In addition, 
for fixed $y>0$ the  $K$-Bessel function
$K_{s}(y)$, which is an entire function in
the $s$ variable, 
is known to have all its zeros on the imaginary axis, 
a result first shown by Polya \cite[p. 308]{Po26}.
(P\'{o}lya's result is stated in terms of  $J$-Bessel
functions, but is identifiable with a  $K$-Bessel function
using well- known identities \cite[(9.6.2), (9.6.3)]{AS}.)
Therefore one can ask: does a similar property hold for
the constant term $a_0(y,s)$, which is now a meromorphic
function of $s$? There is  an interesting answer 
for $y \ge 1$. To state it, recall
 that the {\em modified Riemann hypothesis}
for a function asserts that all its zeros are either on
the line $\Re(s)=\frac{1}{2}$ or on the real axis 
in the interval $0 < x < 1$. 
In 1990  D. Hejhal \cite[Prop. 5.3 (f)]{Hej90},
established that the modified Riemann hypothesis
holds for $a_0(y,s)$ for $y \ge 1$  using the Maass-Selberg relations.
Here we extend this result by determining the occurrence
of real zeros.\\

\begin{theorem}~\label{th13}
For each $y \ge 1$  the constant term of the
Eisenstein series 
\begin{equation}~\label{122}
a_0(y,s) := \zeta^{\ast}(2s)y^s + \zeta^{\ast}(2-2s)y^{1-s}
\end{equation}
is a meromorphic function that satisfies the modified
Riemann hypothesis. There  is a critical value  
\begin{equation}~\label{122a}
y^{\ast}:= 4\pi e^{-\gamma} = 7.055507+
\end{equation}
such that the following hold:

(1) All zeros of $a_0(y,s)$ lie on the critical line for
$1 \le y \le y^{\ast}.$

(2) For $y > y^{\ast}$ there are exactly two zeros off
the critical line. These are real simple zeros 
$\rho_y,1- \rho_y$ with $\frac{1}{2} < \rho_y < 1$. 
The zero $\rho_y$
is a nondecreasing function of $y$, and  $\rho_y \to 1$
as $y \to \infty$.
\end{theorem}

Hejhal~\cite[p. 89]{Hej90}  noted that
for $0< y< 1$ the function $a_0(y,s)$ has complex
zeros off the critical line, with arbitrarily large real part. 
Recently Haseo Ki \cite[Corollary 1]{Ki04b} 
obtained a generalization and strengthening  
of Hejhal's results, remarked on further below. In 
a preprint Ki \cite{Ki04c} shows that all non-real
zeros of $a_0(y,s)$ are simple zeros on the
critical line.

The two real zeros off the critical line 
given in Theorem~\ref{th13} 
seem  of some interest because they behave
like  Landau-Siegel ``exceptional zeros'' as $y \to \infty$. 
The occurrence of real zeros for certain  Epstein zeta
functions was observed in the 1960's
by Bateman and Grosswald \cite{BG64},
and a precise description of the
 ``exceptional zero''  phenomenon in this case
was given by Stark \cite{Sta67}.
Epstein zeta functions are identical  (up to a
nonzero factor) with $E(z,s)$ for fixed $z$. However  
the Riemann hypothesis
fails in general for Epstein zeta functions;
they generically  have zeros off the critical line, and
often  have zeros with $\Re(s) >1$; see Arenstorf
and Brewer \cite{AB93} for some numerical examples. 
We do not address the issue of  proving simplicity of 
zeros of $a_0(y,s)$, but Ki's result \cite{Ki04c}
together with Theorem~\ref{th13} 
gives that for each $y \ge 1$
all zeros of $a_0(y,s)$ are simple, 
except for a multiple zero at $s=\frac{1}{2}$ occurring
for $y= y^{\ast}$. 
As discussed in \S4,
the simple zeros  result is potentially provable
along the lines used for $I(T,s)$, 
but would require additional analysis.

We briefly discuss the proofs.
The proof of Theorem~\ref{th12} rests on 
inequalities proved  on a zero-by-zero basis
in a Hadamard product expansion, 
for a linear combination of two shifted functions.
One may trace methods of this kind back to P\'{o}lya~\cite{Po26},
see also de Bruijn \cite{Bru50}.
Here we make a slight change of hypothesis,
considering only  functions that have both a  functional equation
$F(s) = \pm F(1-s)$ and are real
on the real axis, so have 
the reflection symmetry $F(\bar{s})= \overline{F(s)}$.
In the case of automorphic forms, this 
condition corresponds to self-duality. We formalize
the argument  in Theorem~\ref{th21} in \S2.

The proof of  Theorem~\ref{th13}
follows  Hejhal\cite[Prop. 5.3 (f)]{Hej90} in
deducing  the modified Riemann hypothesis 
using the Maass-Selberg relations. The real zeros
are determined  by inequality estimates.
An alternate proof of the modified Riemann hypothesis here
can be given using an extension of the approach of Theorem~\ref{th11},
which is less elegant and relies on numerical calculations.
This alternate  method 
generalizes to give  information on zeros of the
functions 
\begin{equation}~\label{125a}
H(y,s) := p(s) \zeta^{\ast}(s) y^s + p(1-s) \zeta^{\ast}(2-2s)y^{1-s}
\end{equation}
for $y \ge 1$, provided 
$p(s)$ is a polynomial with real coefficients.
It shows that all but finitely
many of the zeros of $H(y,s)$ lie on the critical line, 
that the zeros off the line are confined to a compact set
independent of $y \ge 1$ and  their number is
uniformly bounded for all $y \ge 1$.
We hope to  treat this method elsewhere.

We now review related work.
In the early 1940's P. R. Taylor, a student of
E. C. Titchmarsh, proved a result similar in form to 
Theorem~\ref{th13} for $y=1$. His  work was
published posthumously \cite{Ta45}. He showed that
$\zeta^{\ast}(s+\frac{1}{2}) - \zeta^{\ast}(s- \frac{1}{2})$
has all its zeros on the critical line. Making the
change of variable $s = 2\tilde{s} - \frac{1}{2}$,
which maps the critical line to itself, this asserts
$F(\tilde{s}) = \zeta^{\ast}(2\tilde{s})- \zeta^{\ast}(2\tilde{s}-1)$, 
satisfies the RH; the 
functional equation yields $\zeta^{\ast}(2\tilde{s}-1)= 
\zeta^{\ast}(2-2\tilde{s})$. In a different direction, the
 functions $I(T, s)$ were  considered by
A. I. Vinogradov and L. Taktajan \cite{VT80}
in 1980, who used them  in an interesting  heuristic argument 
in support of the Riemann hypothesis for $\zeta(s)$.
In 1981 D.  Zagier \cite{Za81a} considered
integrals of Eisenstein series against certain nonnegative
measures supported either at collections of (special) points 
or on collections of (special) closed geodesics. He showed that the 
zero sets of the resulting functions
$F_{\mu}(s)$ contained  the zeros of the Riemann zeta
function, and used them to construct a vector space of
functions carrying a representation of $SL(2, \RR)$ including
principal series representations supported at the zeros
of the Riemann zeta function. 
Recently Haseo Ki \cite{Ki04b}
proved general results 
strengthening and extending those of Hejhal \cite{Hej90},
showing  that finite truncations of the fourier expansion
of Eisenstein series (summed from $-N$ to $N$)
with $z \in \QQ(i)$ and
$\Im (z)$ sufficiently large have all but finitely
many of their zeros simple and lying on the critical line.
In a preprint  Ki \cite{Ki04c} proves that the
constant term $a_0(y,s)$ has all its non-real zeros
simple and on the critical line; his method easily
adapts to determine real zeros,
and could be used to give  independent  
proofs of  Theorems \ref{th11}, \ref{th12} and \ref{th13} above.

The contents of the paper are as follows.
In \S2 we give a theorem  allowing one to
deduce that  all zeros are on a line, and then use it to
prove Theorem~\ref{th12}.
In \S3 we deduce Theorem~\ref{th13}.
In \S4 we establish some results on the density of zeros
on the critical line of $I(T, s)$ and $a_0(y,s)$, and
give some numerical data. We show that for $T \ge 1$ the
zeros of $I(T, s)$ are simple, and those with $\Im(\rho) >0$
have their imaginary parts decrease monotonically as $T$ increases.
In \S5 we make concluding remarks and
raise further topics for investigation. 
In an appendix we describe the interpretation
of L. Weng's zeta functions using semistable lattices. 

H. Ki (private communication)
 observes that there  are now three distinct
approaches to proving Theorem ~\ref{th13}. The first
is an extension of the
proof of Theorem~\ref{th12} above, which is based on a 
variant of the P\'{o}lya approach, and requires some numerical
calculations. The second uses the Maass-Selberg relations
as in Lemma~\ref{le41} below, and is the one given here.
The third is the one taken in \cite{Ki04c}, which 
uses the Hermite-Biehler theorem, and also requires some
numerical calculations; this last approach establishes
simplicity of the zeros.  \\

\noindent{\bf Acknowledgment.} 
The authors thank  Prof. Lin Weng for making them aware 
of mutually overlapping results, leading to this
joint work. They thank  Prof. Dennis Hejhal for informing
them of the work of Prof. Haseo Ki; the preprint
\cite{Ki04c} was done around the same time as
our work. They thank L. Takhtajan for
informing us of his work with A. I. Vinogradov.
They thank the reviewer for helpful comments. 
The first author thanks
David Applegate for exploratory  computations
related to Theorem~\ref{th13}. 

%
%
%
%

\section{Proof of Theorem \ref{th12}.}

The general mechanism used in the proof of
Theorem~\ref{th12} is formalized in the
next result. This result is similar in flavor to a lemma
\footnote{ Both Titchmarsh 
and P\'{o}lya state their results 
in terms  a linear change of
variable $H(w) = F( \frac{1}{2} + iw)$,
asserting that $G(w) = H(w+ic) + H(w-ic)$ 
has real zeros if $H(w)$ 
is real on the real axis and has only real zeros.}
in Titchmarsh\cite[p. 280]{Ti86}, which
traces back to a result of P\'{o}lya ~\cite[Hilffsatz II]{Po26}.
Titchmarsh's lemma assumes (in our notation) that the
function $F(s)$ below has all its zeros on the critical line
(and is real there)
which we relax to (ii) by assuming
extra symmetries in  (i). \\

\begin{theorem}~\label{th21}
Let $F(s)$ be an entire function of genus zero or one,
that has the following properties.

(i) $F(s)$ is real on the real axis, and satisfies a functional
equation of form 
\begin{equation}~\label{201}
F(s) = \pm F(1-s),
\end{equation}
for some choice of sign.

(ii) There exists $a > 0$ such
that all  zeros of $F(s)$ lie in the vertical strip
\begin{equation}~\label{202}
 \left| \Re(s)- \frac{1}{2} \right| < a.
\end{equation}

Then for any real $c \ge a$,
\begin{equation}~\label{203}
\left| \frac{F(s+ c)}{F(s-c)} \right| > 1~~\mbox{~~if~~} 
\Re(s) > \frac{1}{2},
\end{equation}
and
\begin{equation}~\label{204}
\left| \frac{F(s+ c)}{F(s-c)} \right| < 1~~\mbox{~~if~~} 
\Re(s) < \frac{1}{2}.
\end{equation}
In particular, for any $0 \le \theta <  2 \pi$ the function
\begin{equation}~\label{205}
G_{\theta}(s) := F(s+c) + e^{i \theta} F(s-c)
\end{equation} 
has all its zeros falling on the line $\Re(s) = \frac{1}{2}$.
\end{theorem}

\paragraph{{\bf Proof.}}
The genus one assumption is equivalent to the assertion that
the Hadamard product factorization
\begin{equation}~\label{206}
F(s) = e^{A+Bs} s^R \prod_{\rho}( 1 - \frac{s}{\rho}) e^{\frac{s}{\rho}}
\end{equation}
converges absolutely and uniformly on compact
subsets of $\CC$, see Titchmarsh \cite[Sec. 8.23, 8.24]{Ti39}.
This  assumption is also equivalent to the bound
$\sum_{\rho}  \frac{1}{|\rho|^2} < \infty.$
 Hypothesis (i) implies  symmetries of the zeros around both
the real axis and the line $\Re(s) = \frac{1}{2}$; i.e. under
both $\rho \mapsto 1- \rho$ and $\rho \mapsto \bar{\rho}$.
It follows that  set of zeros $\rho= \beta + i \gamma$,
counted with multiplicity, can be partitioned  into  blocks $B(\rho)$
comprising  
$\{ \rho, 1-\rho, \bar{\rho}, 1 - \bar{\rho} \}$ if $\beta \ne \frac{1}{2}$;
$\{ \rho, 1- \rho \}$ if $\beta= \frac{1}{2}$ and $\gamma \ne 0$;
and 
$\{\rho\}$ if $\rho= \frac{1}{2}$. Each block is labelled with the
unique zero in it having $\beta \le \frac{1}{2}$ and $\gamma \ge 0$.
We next assert, using hypothesis (ii),  that 
the modified Hadamard product obtained by 
grouping over blocks, can have
the convergence factors $e^{\frac{s}{\rho}}$ removed, i.e. 
\begin{equation}~\label{207}
F(s) = e^{A+ B' s} \prod_{B(\rho)} 
\left( \prod_{\rho \in B(\rho)}( 1 - \frac{s}{\rho})\right), 
\end{equation}
where the outer product on the right side 
converges absolutely and uniformly on compact
subsets of $\CC$. This assertion holds because the block convergence
factors $e^{c(B(\rho))s}$ are given by 
\begin{eqnarray*}
c(B(\rho)) &= &  \frac{\beta}{|\rho|^2} + \frac{1-\beta}{|1 - \rho|^2}
\mbox{~~if~~} \beta \ne \frac{1}{2}; \\ 
c(B(\rho)) &=& \frac{1}{|\rho|^2} \mbox{~~if~~} \beta = \frac{1}{2},
\gamma \ne 0; \\
c(B(\rho))& = & \frac{\frac{1}{2}}{|\rho|^2}=2  
\mbox{~~if~~}  \rho= \frac{1}{2}.
\end{eqnarray*} 
Hypothesis (ii) gives $- a < \beta-\frac{1}{2} < a$ hence 
$$
\sum_{B(\rho)} |c(B(\rho))| \le 
(1+2a) \left( \sum_{\rho} \frac{1}{|\rho|^2} \right) < \infty.
$$
Thus the convergence factors can be pulled out of the product,
 yielding \eqn{207}, with $B' = B + \sum_{B(\rho)}c(B(\rho)).$

Using the functional equation \eqn{201} we can further
infer that the constant $B'=0$ in \eqn{207}, so that 
\begin{equation}~\label{207a}
F(s) = e^{A} \prod_{B(\rho)} 
\left( \prod_{\rho \in B(\rho)}( 1 - \frac{s}{\rho})\right).
\end{equation}
Indeed the  change of variable $s \mapsto 1-s$ permutes the
factors in each block $B(\rho)$, with a possible sign change
for $\rho= \frac{1}{2}$,  
so it must be that
$e^{A+ B's}= \pm e^{A+B'(1-s)}$, which forces $B'=0$.

To establish \eqn{203} and \eqn{204} we can now proceed
block by block in \eqn{207a}, using the factorization 
\begin{equation}~\label{208}
\left| \frac{F(s+ c)}{F(s-c)} \right| = 
\prod_{B(\rho)} \left( \prod_{\rho \in B(\rho)}
\left| \frac{1- \frac{s+c}{\rho} }{1- \frac{s-c}{\rho} } \right|\, \right).
\end{equation}
In a single block we can clear denominators to obtain 
$$
\prod_{\rho \in B(\rho)}
\left| \frac{1- \frac{s+c}{\rho} } {1- \frac{s-c}{\rho} } \right| = 
\prod_{\rho \in B(\rho)} \left| \frac{s+c -\rho}{s-c -\rho} \right|.
$$
The main point is now to compare the  term in the numerator
with $\rho$ against the term in the denominator
with $\rho':= \overline{ 1 - \rho}= 1-\bar{\rho}$. 
We assert that 
\begin{equation}~\label{209}
\left| \frac{s+c - \rho} {s-c - (1-\bar{\rho})  } \right|^2 > 1
\mbox{~~if~~~} \Re(s) > \frac{1}{2},
\end{equation} 
and
\begin{equation}~\label{210}
\left| \frac{s+c - \rho} {s- c - (1-\bar{\rho} ) } \right|^2 < 1
\mbox{~~if~~~} \Re(s) < \frac{1}{2}.
\end{equation}
If \eqn{209} is shown, then we may conclude
for $\Re(s) > \frac{1}{2}$  that the
absolute value of the product over terms in each block
on the right in \eqn{208} exceeds $1$, and \eqn{203}
follows. Similarly \eqn{210} implies that for
 $\Re(s) < \frac{1}{2}$  the product of terms over
each block is smaller than $1$, and \eqn{204}
follows.

It remains to show \eqn{209} and \eqn{210}.
Writing $s= \sigma + it$, we have
$$
\left| \frac{s+c - \rho}{s-c - (1-\bar{\rho})} \right|^2=
\frac{(\sigma + c- \beta)^2 + (t-\gamma)^2}
{(\sigma -c - 1 + \beta)^2 + (t-\gamma)^2}
$$
Now \eqn{209} reduces to the assertion that
\begin{equation}~\label{211}
(\sigma + c - \beta)^2 > (\sigma -c - 1 + \beta)^2 
\mbox{~~when~~} \Re(s) > \frac{1}{2}.
\end{equation}
To show this we note that $\Re(s) > \frac{1}{2}$ gives
$$
\sigma+ c - \beta > \frac{1}{2}+ a  -\beta > 0, 
$$
whence  \eqn{211} makes the two assertions
$$
\sigma + c- \beta > \sigma-c-1 + \beta,
$$
$$
\sigma + c- \beta > -(\sigma-c-1 + \beta).
$$
The second of these  asserts that $\sigma > \frac{1}{2}$,
while the first asserts that  $2c > 2(\beta - \frac{1}{2})$,
which holds since $c \ge a > \beta - \frac{1}{2}.$
Thus \eqn{211} holds, whence \eqn{209} holds. 
  
 A similar argument is used to establish  
\eqn{210}. It reduces to the
assertion that
\begin{equation}~\label{212}
(\sigma + c - \beta)^2 < (\sigma -c - 1 + \beta)^2 
\mbox{~~when~~} \Re(s) < \frac{1}{2}.
\end{equation}
We have
$$ 
-(\sigma -c - 1 + \beta) \ge \frac{1}{2} + a - \beta > 0,
$$
so that  \eqn{212} is equivalent to the two assertions
$$
-\sigma + c + 1 - \beta > -(\sigma + c - \beta),
$$
$$
-\sigma + c + 1 - \beta > \sigma + c - \beta.
$$
The second of these is equivalent to
$\sigma < 1/2$ and the first to
$c + 1/2 - \beta > 0$, which holds. 

The conclusion that $G_{\theta}(s)$ has all its zeros on the
line $\Re(s) = \frac{1}{2}$ follows, because
the two terms on the right side of \eqn{205} have differing
absolute values off this line. 
$~~~\Box$\\

\paragraph{\bf Remarks.} (1) In the special case  when 
$e^{i\theta} = \pm 1$,  the function  
$G_{\theta}(s)$ has a  functional equation
\begin{equation}~\label{213}
G_{\theta}(s) = \pm G_{\theta}(1-s)
\end{equation}
inherited from the functional equation  of $F(s)$.
This need not hold for other values of $\theta$. \\

(2) The proof of Theorem~\ref{th21} has two main
steps: 
first, establishing the existence of a
 modified Hadamard product  factorization grouping zeros
into finite blocks, with the convergence factors
dropped, and second, a local inequality argument that
applies to each block separately.
The local inequality argument
used in the proof paired zeros $\rho$ and $1 - \bar{\rho}$,
a condition that requires only 
that $F(s)$ have constant modulus on the
critical line, and does not require symmetry
around the real axis. This condition is still met when 
hypothesis (i) is  relaxed  to: \\

{\em (i') $F(s)$ satisfies 
a functional equation of form
\begin{equation}~\label{213b}
F(s) = e^{i \alpha} F(1-s)
\end{equation}
for some $0 \le \alpha < 2\pi$.}  \\

\noindent
The hypothesis (i') allows general $L$-functions.
In the argument above hypothesis (i) was used to get the
modified Hadamard product factorization. It is 
possible to prove an alternative version 
of Theorem~\ref{th21} assuming (i') and (ii),
provided  some extra hypothesis  (iii) is added
restricting the locations of zeros, that guarantees the
existence of  a modified Hadamard product 
\eqn{207} that  converges absolutely and uniformly
on all compact subsets of the plane. When (i') holds  
the zeros no longer need to be symmetric about the
real axis, so  blocks must be chosen differently
to get convergence; such blocks remain invariant 
under the map $\rho \mapsto 1- \bar{\rho}.$
The functional equation (i') then implies only that
$\Re(B')=0$, but this is sufficient to obtain the result.

It can be shown that  automorphic $L$-functions
(principal $L$-functions for $GL(n)$ ) do possess such modified
Hadamard product expansions, using 
asymptotic formulas for zeros 
given in \cite[Theorem 2.1 (4)]{La04b}. \\


\paragraph{{\bf Proof of Theorem~\ref{th12}.}}

Let $\xi(s)$ be the Riemann $\xi$-function
\begin{equation}~\label{220-1}
\xi(s):=\frac{1}{2}s(s-1)\zeta^\ast(s),
\end{equation}
which is an entire function. 
Theorem~\ref{th21} applies to 
$F(s) = \xi(2s-\frac{1}{2})$, which satisfies
the functional equation
$F(s) = F(1-s)$, is real on the real axis,
and the zeros
of $F(s)$ are confined to $\frac{1}{4} < \Re(s) < \frac{3}{4}$,
so we can take  $a= \frac{1}{4}.$

Now we apply Theorem \ref{th21} with  $c= \frac{1}{4}$,
which for the function 
$$G(s) := F(s+ \frac{1}{4}) + F(s-\frac{1}{4}) = \xi(2s) + \xi(2s-1)$$
yields, for $\Re(s) > \frac{1}{2}$,  that
\begin{equation}~\label{220}
\left| \frac{F(s+ \frac{1}{4})}{F(s-\frac{1}{4})} \right| =
\left| \frac{\xi(2s)}{\xi(2s-1)} \right| > 1 
\mbox{~~for~~} \Re(s) > \frac{1}{2}.
\end{equation}

Recall now that
$$
I(T, s) := -\frac{\zeta^{\ast}(2s)}{s-1} T^{s-1} + 
\frac{\zeta^{\ast}(2-2s)}{s} T^{-s}
= -\frac{\zeta^{\ast}(2s)}{s-1} T^{s-1} + 
\frac{\zeta^{\ast}(2s-1)}{s} T^{-s},
$$
using the functional equation $\xi(s) = \xi(1-s)$.
For fixed real $T>0$ the  function $I(T, s)$
has simple poles at $s=0, 1$ and is analytic elsewhere,
and satisfies the functional equation $I(T, s) = I(T, 1-s)$.
At $s=1$ it  has residue $ c_1(T)=-\frac{1}{2}(\frac{\pi}{3} - \frac{1}{T})$,
  and at $s=0$ residue $-c_1(T)$;
both terms in 
the expression  \eqn{120} for $I(T,s)$
contribute to these residues.
To prove the theorem it suffices to study zeros of the 
entire function 
\begin{equation}~\label{222}
H(T,s) := \frac{1}{4} (2s)(2s-1)(2s-2) I(T, s)
= -\xi(2s) T^{s-1} + \xi(2s-1)T^{-s}.
\end{equation}
This function  has a zero at $s= \frac{1}{2}$ and satisfies the
functional equation
\begin{equation}~\label{223}
H(T,s) = - H(T, 1-s). 
\end{equation}
Applying the result \eqn{220}, we have for $\Re(s) > \frac{1}{2}$
and $T \ge 1$ that 
$$
\left| \frac{ -\xi(2s) T^{s-1} } {\xi(2s-1)T^{-s}} \right| =
\left| \frac{ \xi(2s)}{\xi(2s-1)} \right| T^{ 2 \sigma -1} \ge 
\left| \frac{ \xi(2s)}{\xi(2s-1)} \right| > 1.
$$
We conclude for $T \ge 1$ that $H(T, s) \ne 0$ when
$\Re(s) > \frac{1}{2}$, and  the functional equation \eqn{223}
then yields $H(T, s) \ne 0$ when $\Re(s) < \frac{1}{2}.$ 
Thus for $T \ge 1$ all zeros of $H(T, s)$ must have
$\Re(s) = \frac{1}{2}.$
$~~~\Box$



\section{Proof of Theorem~\ref{th13}.}

The simplest entire function associated to
$a_0(y,s)$ is
\begin{equation}~\label{400a}
G(y, s) := (2s)(2s-2) a_0(y,s),
\end{equation}
which behaves similarly
to the Riemann $\xi$-function, satisfying
the functional equation  $G(y,s) = G(y, 1-s)$, 
being real on the real axis and on the line $\Re(s) = \frac{1}{2}$.
It also has $G(y,\frac{1}{2}) = (\log 4 \pi - \gamma-\log y)\sqrt{y}$,
where $\gamma$ is Euler's constant.
However to  establish  Theorem~\ref{th13} it proves useful
to study instead the entire function  
\begin{equation}~\label{400}
H(y, s) := \frac{1}{2}(s- \frac{1}{2}) G(y,s)= (s-1)\xi(2s)y^s +
s\, \xi(2s-1)y^{1-s}.
\end{equation}
which adds an extra zero at $s= \frac{1}{2}$.
The function $H(y,s)$ 
satisfies the functional equation
$H(y,s) = -H(y,1-s)$, but has
the advantage that  both terms on the right side 
of \eqn{400} are entire functions.

First we will show using the Maass-Selberg relations
( Hejhal \cite[Prop. 5.3(f) ]{Hej90}) that 
all zeros of $a_0(y,s)$ lie on the critical line for any $y \geq 1$ 
except for real zeros.
Second we will determine  the location of real zeros of $a_0(y,s)$. \\




\begin{lemma}[Hejhal]~\label{le41}
For any $y \geq 1$, 
all zeros of $a_0(y,s)$ lie on the critical line except for real zeros.  
\end{lemma}

\paragraph{\bf Proof.} 

Let 
\begin{equation}~\label{401}
E_T^\ast(z,s):=
\begin{cases}
E^\ast(z,s) -a_0(y,s)
& \text{ if $z \in D^T$ }, \\
E^\ast(z,s)
& \text{ if $z \in D - D^T$ }
\end{cases}
\end{equation}

where $T \geq 1$. 
The Maass-Selberg relation is stated as (cf.\cite[pp.154--155]{Hej83})
\begin{align}~\label{402}
\aligned
(s-\bar{s})(1-s-\bar{s}) \int_{D} 
 & |E_T^\ast(z,s)|^2 \frac{dxdy}{y^2} \\
= &  a_0(T,s) (\bar{s} \zeta^\ast(2\bar{s}) T^{\bar{s}-1} + 
(1-\bar{s})\zeta^\ast(2\bar{s}-1)T^{-\bar{s}}) \\
& - a_0(T,\bar{s}) (s \zeta^\ast(2s) T^{s-1} + (1-s)\zeta^\ast(2s-1)T^{-s}). 
\endaligned
\end{align}
If $a_0(T,s)=0$ then,  by the reflection
principle, $a_0(T,\bar{s})=0$ as well, so that the right side of \eqn{402}
is zero. But the left side 
of \eqn{402} is non-zero whenever $\Re(s) \ne 0$ 
and $\Re(s) \not= 1/2$ both hold. 
$~~~\Box$ \\

Lemma~\ref{le41} shows that we need only to determine 
the locations of the zeros of $H(y,s)=s(s-1)(2s-1)a_0(y,s)$ 
on the real axis.  \\


\begin{lemma}~\label{le42} 
For any fixed $y>0$, the constant term $a_0(y,\sigma)$ has no zero 
outside of the open interval $(0,1)$ as a function of $\sigma$ on $\RR$.   
\end{lemma}

\paragraph{\bf Proof.} 
For real $\sigma>1$ we have $E^\ast(z,\sigma)>0$ from the series 
representation $\eqn{101}$. 
If $a_0(y,\sigma_0)=0$ for some $\sigma_0>1$, then the integral 
$\int_0^1 E^{\ast}(x+iy,\sigma_0) dx$ is equal to $0$. 
This is a contradiction. Thus $a_0(y,\sigma_0)\ne 0$ for $\sigma > 1$ and
the same holds for $\sigma < 0$ using the functional equation.
$~~~\Box$ \\ 

Using  Lemma~\ref{le42} and the functional equation it sufficies to deal 
with $H(y,s)$ for real $s= \sigma$ in the  interval $(\frac {1}{2},1)$. 
To prove the theorem it
suffices to show  that $H(y,\sigma)\not=0$ on $(\frac{1}{2},1)$ 
for any fixed $ 1 < y \leq y^\ast =4\pi e^{-\gamma}$ 
and that there is only one zero of $H(y,\sigma)$ in $(\frac{1}{2},1)$ 
for any fixed $y > y^\ast$. 

We can rewrite the condition  $H(y,\sigma)=0$ as 
\begin{equation}~\label{403}
\frac{1-\sigma}{\sigma} \frac{\xi(2\sigma)}{\xi(2\sigma-1)} y^{2\sigma-1}=1.
\end{equation}
Here we consider the function
\begin{equation}~\label{40a3}
F(y,\sigma):= y^{2 \sigma-1 }f(\sigma) \quad \text{with} \quad 
f(\sigma) := \frac{1-\sigma}{\sigma} \frac{\xi(2\sigma)}{\xi(2\sigma-1)}.
\end{equation}
From the product formula $\eqn{207}$ we find that $f(\sigma)>0$
for $0 <\sigma <1$. 
In fact we see that
\begin{equation}~\label{405}
f(\sigma)
=\frac{1-\sigma}{\sigma}
\prod_{B(\rho)}
\left(\prod_{{\rho \in B(\rho)}\atop{\gamma>0}} 
\frac{(2\sigma-\beta)^2+\gamma^2}{(2\sigma-2+\beta)^2+\gamma^2}\right) >0.
\end{equation}

Now let 
\begin{equation}~\label{404}
F^\prime (y,\sigma) :=\frac{d}{d \sigma} F(y,\sigma).
\end{equation}
Then we have
\begin{equation}~\label{404a}
F^\prime (y,\sigma) =2y^{2\sigma-1}f(\sigma) \, 
\left( \log y +\frac{1}{2} \frac{f^\prime(\sigma)}{f(\sigma)} \right).
\end{equation}
Now \eqn{405} implies
that for $y>0$ the function $F^\prime (y,\frac{1}{2})$ is a 
strictly increasing function of $y$. 
Here we establish the following lemma. \\


\begin{lemma}~\label{le43} 
\begin{enumerate} 
\item[(a)] For $y > 0$, $F^\prime (y,\frac{1}{2})=0$ at the unique value 
$y^\ast=4\pi e^{-\gamma}=7.055507+$. 

\item[(b)] The function $-f^\prime(\sigma)/f(\sigma)$ is a 
strictly increasing function of $\sigma$ on $(\frac{1}{2},1)$.  
\end{enumerate}
\end{lemma}

\paragraph{\bf Proof.} 
First we prove (a). 
Because $f(\sigma)>0$ for $0< \sigma <1$, 
$F^\prime(y,\frac{1}{2})=0$ implies that 
$2 \log y=-f^\prime(\frac{1}{2})/f(\frac{1}{2})$. 

We have
\begin{align}
-\frac{f^\prime(\frac{1}{2})}{f(\frac{1}{2})}
= 4 -2\frac{\xi^\prime(1)}{\xi(1)} +2\frac{\xi^\prime(0)}{\xi(0)} 
= 4(1+\frac{\xi^\prime(0)}{\xi(0)}),
\end{align}
using $\xi(0) = \xi(1)$ and $\xi^\prime(0) = -\xi^\prime(1)$.
Therefore the unique value $y^\ast$ where $F^\prime(y,\frac{1}{2})=0$ 
is given by 
$$
\log y^\ast = 2(1+\frac{\xi^\prime(0)}{\xi(0)}).
$$ 
We recall the fact that 
$$
\frac{\xi^\prime(0)}{\xi(0)}=\frac{1}{2}\gamma -1 + \frac{1}{2}\log 4\pi 
= -0.0230957+,
$$
where $\gamma=0.57721+$ is Euler's constant,
see Davenport ~\cite[pp. 80--82]{Da80}. Hence $y^\ast=4 \pi e^\gamma$.

Next we prove (b). We have
\begin{equation}~\label{406}
\left( - \frac{f^\prime(\sigma)}{f(\sigma)} \right)^\prime
=\frac{2\sigma-1}{\sigma^2(1-\sigma)^2} 
+2 \left( - \frac{\xi^\prime}{\xi}(2\sigma) + 
\frac{\xi^\prime}{\xi}(2\sigma-1) \right)^\prime.
\end{equation}

The first term in the right hand side in $\eqn{406}$ is
 positive for $\frac{1}{2}< \sigma <1$, so it suffices 
to show that the second term in the right hand side 
in $\eqn{406}$ is also positive for $\frac{1}{2}< \sigma <1$. 

By using the product formula $\eqn{207}$ we obtain
\begin{align}
\aligned
- \frac{\xi^\prime}{\xi}(2\sigma) + \frac{\xi^\prime}{\xi}(2\sigma-1)
& = \sum_{B(\rho)} \sum_{\rho \in B(\rho)} 
\left( -\frac{1}{2\sigma-\rho}+\frac{1}{2\sigma-1-(1-\rho^\prime)} \right) \\
& = \sum_{B(\rho)} \sum_{{\rho \in B(\rho)}\atop{\gamma >0}} 
\left( -\frac{2(2\sigma-\beta)}{(2\sigma-\beta)^2+\gamma^2}
+\frac{2(2\sigma-2+\beta)}{(2\sigma-2+\beta)^2+\gamma^2)} \right), 
\label{406a}
\endaligned
\end{align}
where $\rho=\beta+i\gamma$. 

We show that the derivatives of each term in the second line 
of $\eqref{406a}$ 
are positive for $\frac{1}{2} < \sigma <1$. 
Take
\begin{equation}~\label{407}
g(\sigma):=\frac{d}{d \sigma} 
\left( 
-\frac{2(2\sigma -\beta)}{(2\sigma -\beta)^2 +\gamma^2}
+\frac{2(2\sigma -2 +\beta)}{(2\sigma -2 +\beta)^2 +\gamma^2}
\right). 
\end{equation}

Then we have
\begin{align*} \aligned
g(\sigma) & = -\frac{ 4\gamma^2 -4(2\sigma-\beta)^2 }
{ \left((2\sigma-\beta)^2 +\gamma^2\right)^2 }
+\frac{ 4\gamma^2 +4(2\sigma-2+\beta)^2}
{ \left((2\sigma-2+\beta)^2+\gamma^2\right)^2} \\
& = \frac{ 4(2\sigma-1) (1 -\beta) 
[\, 4\gamma^2 \{ (2\sigma -\beta)^2+(2\sigma-2 +\beta)^2 \} 
+12\gamma^4 -4(2\sigma -\beta)^2(2\sigma-2 +\beta)^2 \,] }
{ \left((2\sigma-\beta)^2 +\gamma^2 \right)^2 \, 
\left(2\sigma-2+\beta)^2+\gamma^2\right)^2 }.
\endaligned  \end{align*}
Because $0<\beta<1$ 
$$
4(2\sigma-\beta)^2 (2\sigma-2+\beta)^2 <16
$$
for $\frac{1}{2}<\sigma<1$ and we know that $\gamma>14$. 
Hence $g(\sigma)>0$ for $\frac{1}{2}<\sigma<1$. 
This implies $(-f^\prime(\sigma)/f(\sigma))^\prime>0$ 
for $\frac{1}{2}<\sigma<1$, which proves the
lemma.  $~~~\Box$\\

To complete the proof of Theorem~\ref{th13}, 
from Lemma~\ref{le43} (b), we have $F^\prime(y,\sigma)=0$ 
for at most 
one $\sigma \in (\frac{1}{2},1)$ for any fixed $y \geq 1$. 
Now suppose  $1 \leq y \le y^\ast$. We then have 
\begin{equation}
\log y \le \log y^\ast 
= -\frac{1}{2}\frac{f^\prime(\frac{1}{2})}{f(\frac{1}{2})} 
< -\frac{1}{2}\frac{f^\prime(\sigma)}{f(\sigma)} 
\end{equation}
for $\frac{1}{2} < \sigma <1$. Therefore
\begin{equation}
\log y + \frac{1}{2}\frac{f^\prime(\sigma)}{f(\sigma)} < 0
\end{equation}
so that $F^\prime(y,\sigma)<0$ on $(\frac{1}{2},1)$. 
Since $F(y,\frac{1}{2})=1$, $F(y,\sigma)\not=1$ for any 
$\frac{1}{2} < \sigma <1$, and this implies that  
$H(y,\sigma)\not=0$ on  $(\frac{1}{2},1)$. 

Next suppose  $y > y^\ast$. Then there is
a unique $\sigma_0$ with $\frac{1}{2} < \sigma_0 < 1$ such that 
\begin{equation}
\log y^\ast < \log y 
= -\frac{1}{2}\frac{f^\prime(\sigma_0)}{f(\sigma_0)}, 
\end{equation}
because $-f^\prime(\sigma)/f(\sigma) \to +\infty$ 
monotonically as $\sigma \to 1$. 
Further
$$  -\frac{1}{2}\frac{f^\prime(\sigma_1)}{f(\sigma_1)} 
< -\frac{1}{2}\frac{f^\prime(\sigma_0)}{f(\sigma_0)} 
< -\frac{1}{2}\frac{f^\prime(\sigma_2)}{f(\sigma_2)}~\mbox{for} ~
\frac{1}{2} \leq \sigma_1 < \sigma_0 < \sigma_2 <1,
$$
hence 
$F^\prime(y,\sigma_1) >0$ for $\frac{1}{2} \leq \sigma_1 <\sigma_0$
and 
$F^\prime(y,\sigma_2) <0$ for $\sigma_1< \sigma_2 <1$.
Since $F(y,\frac{1}{2})=1$ and $F(y,1)=0$, 
these imply there is a unique value $\sigma_y$ in $(\frac{1}{2}, 1)$
such that $F(y,\sigma_y)=1$, and this value has $\sigma_0 <\sigma_y < 1$. 
This is exactly the condition for $H(y, \sigma)=0$ so
we conclude that the unique value where this occurs in $(\frac{1}{2}, 1)$
is $\sigma= \sigma_0$. We also  find 
~\footnote{ Asymptotically we see that 
$\sigma_0=1+O(\frac{1}{\log y})$. 
This suggests that $\sigma_y=1+O(\frac{1}{\log y})$. } 
that $\sigma_0 \to 1$ as $y \to +\infty$, 
which implies that  $\sigma_y \to 1$ as $y \to +\infty$. 
This completes  the proof of Theorem~\ref{th13}. $~~~~\Box$ \\

%
%
%
%
\section{Distribution of Zeros}

We have shown in \S2 that the Riemann hypothesis
holds for $I(T,s)$ for each fixed $T \ge 1$ and
in \S3 that the modified Riemann hypothesis holds
for $a_0(y,s)$ for each $y \ge 1$. 
Here we consider how the zeros behave as the
parameter $T$ or $y$ is varied, and also consider the issue of 
simplicity of zeros.
Note that the zeros of these functions vary continuously in $T$ as the 
parameter $T$  is varied, and vary  analytically in $T$ as long as they
are simple zeros.

In what follows we let
$N(f, U)$ count the number of zeros of the function $f(s)$
having $|\Im(s) | \le U$.

\begin{theorem}~\label{Nth51}
(1) For each fixed $T \ge 1$ the function $I(T,s)$ has simple zeros.
The number of zeros of $I(T,s)$ with $|\Im(s)| \le U$
satisfies 
\begin{equation}~\label{N500}
N(I(T,s); U) = N( \xi(2s); U) + \frac{2}{\pi} (\log T) U + O(\log U).
\end{equation}

(2) As $T\ge 1$ increases, each  zero $\rho$
 of $I(T,s)$ with $\Im(\rho) >0$
has imaginary part that is a
strictly decreasing function of $T$.
\end{theorem}

\paragraph{\bf Remark.} Standard estimates from the theory of 
the Riemann zeta function show that
\begin{equation}~\label{N501}
 N(\xi(2s), U) = \frac{2}{\pi}U \log U -
\frac{2}{\pi}\left(\log \pi +1\right) U + O (\log U).
\end{equation}
see for example \cite[p. 98]{Da80}.

\paragraph{\bf Proof.}
We study the variation in argument of the 
 entire function
$$
H(T,s) :=\frac{1}{4}(2s)(2s-1)(2s-2) I(T,s) = 
-\xi(2s) T^{s-1} +\xi(2-2s)T^{-s},
$$
on the critical line $s= \frac{1}{2} + it$. The function $H(T,s)$
has the same zeros as $I(T,s)$, with multiplicity,
 except for an extra zero at
$s=\frac{1}{2}$ (where $I(T, \frac{1}{2}) \ne 0$)
 and by Theorem~\ref{th12} all
its zeros lie on $\Re(s)= \frac{1}{2}.$
It is pure imaginary-valued on the critical line,
with 
$$
H(T,\frac{1}{2} + it ) = -\xi(1+ 2it) T^{-\frac{1}{2} + it}
+ \xi(1-2it) T^{-\frac{1}{2} - it}.
$$
Now set $\xi(1+it)= R(t) e^{i \theta(t)} $ with $R(t) = |\xi(1+it)|$,
where the argument $\theta(t)$ is measured continuously from $\theta(0)=0$
on the real axis. 
Then we have
$$
H(T,\frac{1}{2} + it )= -\frac{2i}{\sqrt{T}}R(2t)\sin(\theta(2t) + t\log T ),
$$
We  detect zeros exactly when
\begin{equation}~\label{N503}
f(T,t) := \theta(2t) + t \log T \equiv 0 ~(\bmod~\pi),
\end{equation}
but a priori have no information on their multiplicity.

(1) To show the  zeros detected by \eqn{N503} are simple zeros of 
$H(T,s)$, we use results from Lagarias \cite{La04}.
We treat the cases $T=1$ and $T> 1$ separately. 
For $T=1$ we use a special case 
of Theorem 2.1 of \cite{La04}, which  asserts that
$$
\xi_{1/2, \pi/2}(s) := - \frac{1}{2i}
\left( \xi(s+ \frac{1}{2}) - \xi(s- \frac{1}{2})\right)
$$ 
has simple zeros which
all lie on the critical line $\Re(s) = \frac{1}{2}$. 
Using the functional equation
we have
$$
 \xi_{1/2, \pi/2}(2s-\frac{1}{2}) = - \frac{1}{2i}
\left( \xi(2s) - \xi(2s- 1)\right) =- \frac{1}{2i}
\left( \xi(2s) - \xi(2-2s)\right) = - \frac{1}{2i}H(1, s)
$$
It follows that all  the zeros of $H(1,s)$ are simple and lie on the
critical line.

For $T >1$ it is sufficient to show that 
$\frac{d}{dt}f(t,T)>0$ for all $t \in \RR$. We have
$\frac{d}{dt}f(t,T) = \frac{d}{dt}f(t,1) + \log T$,
and since $\log T >0$ 
it suffices to  establish $\frac{d}{dt}f(t,1) \ge 0$
for all $t$.
We now apply Lemma 2.2 of \cite{La04} to $E(s)= \xi(s + \frac{1}{2})$,
whose hypotheses are met
by Lemma 2.1 of that paper. Now equation (2.10) of that
paper yields $\frac{d}{dt} \theta(t) \ge 0$ for all
$T$, as required. 

It follows that $H(T, s)$ has simple
zeros for each $T \ge 1$, hence so does $I(T,s)$.
For fixed $T$ we  now may number the
zeros of $H(T,s)$ 
with $\Im(\rho)>0$ in order of  increasing imaginary part as
$0< \gamma_1(T) < \gamma_2(T) < ...$; and those with $\Im(\rho)<0$
are numbered by  $\gamma_{-n}(T) = - \gamma_n(T)$; here $\gamma_{0}(T) = 0$
is the zero of $H(T,s)$ at $s= \frac{1}{2}.$ 

To count zeros, we have from \eqn{N503} that 
\begin{equation}~\label{N504a}
N(I(T,s); U) = \frac{2}{\pi} f(T,U) + O(1) =
\frac{2}{\pi} f(1,U) + \frac{2}{\pi} (\log T) U + O(1).
\end{equation}
The function $f(1,U) = \theta(2U)$ measures 
the change in argument of $\xi(s)$, from $s=1$ to $s= 1+ 2iU$. 
Following  Davenport \cite[Chap. 15]{Da80},
this argument change differs from that of $\xi(\frac{1}{2} + 2iU)$ by
$O(\log U)$, and the
argument change to this point  is $\frac{1}{4} (2\pi) N(\xi(2s), U) +O(1)$,
since the value $(2\pi) N(\xi(2s), U) $ is a contribution
of four terms of this type. We conclude that
$\frac{2}{\pi}f(1,U) =N(\xi(2s),U) + O(\log U)$.
Substituting this in \eqn{N504a} yields \eqn{N500}.

(2) This fact follows directly from \eqn{N503}. That is,
$f(T,t)$ is strictly increasing in both $t$ and $T$, so
increasing $T$ increases the rate of turning of the function,
so the $n$-th time the value $0 (\bmod~\pi)$ is reached is
a strictly decreasing function of $T$.
$~~~~\Box$ \\

The first few zeros of $I(T, s)$ for
$T=1$ and $T=y^{\ast}=  4\pi e^{-\gamma} = 7.055507+ $
are given in Table \ref{zeros1} below. (The value $T=y^{\ast}$
was chosen for comparison with the zeros of $a_0(y,s)$ 
zeros at this point.)

\begin{table}[htbp]
\center
\begin{tabular}{|r|c|c|}
\hline
 & $T=1$ & $T=y^{\ast}$ \\
\hline
1 & 7.769080112 & 1.570199673 \\
2 & 11.01900402 &  3.136172650 \\
3 & 13.11079833 & 4.688303082  \\
4 & 15.58052582 & 6.172131737 \\
5 & 17.07367093 & 7.073883755 \\
6 & 19.21539818 & 7.990000858 \\
7 & 20.83659268 & 9.408931700  \\
8 & 22.24754162 & 10.42818054  \\
9 & 24.26973459 & 11.18590514  \\
10 & 25.39649716 &12.29085328  \\
11 & 26.95610030 &12.94293137  \\
12 & 28.59571466 &14.21495758 \\
13 & 29.93119639 &15.15516152  \\
14 & 31.03561085 &15.87628765  \\
15 & 32.83737170 & 16.56289948  \\
\hline
\end{tabular}
\caption{Zeros of constant term $I(T, s)$ on critical line}
\label{zeros1}
\end{table}

For the  constant term $a_0(y,s)$
similar properties of its zeros should hold, with minor
modifications. One can show by similar arguments that 
$$
N( a_0(y,s), U) = N( \xi(2s), U) + \frac{4}{\pi} \log y U + O(1).
$$
The motion of the 
zeros as $y\ge 1$ is slightly more complicated than
for $I(y,s)$, 
since a multiple zero occurs at $y=y^{\ast}$, 
and two zeros eventually migrate
to the real axis. It appears
that, aside from these two zeros,  all zeros are simple and 
their ordinates monotonically decrease as
the parameter $y$ increases. To prove
this rigorously one have to  consider instead of 
the argument of $\xi(2s) y^{s-1}$ the
argument of $(2s-2) \xi(2s)y^s$ on the
critical line,  and the argument
of $-1+2it$ turns in the wrong direction. 
Numerical and analytic estimates would be needed to control
this  effect.
Presumably the contribution of $-1+2it$ and its 
conjugate to the argument accounts for
the escape of the two zeros to the real axis.
The first few complex zeros of $a_0(y,s)$ are given 
the values $y=1$ and $y= y^{\ast}$ in Table \ref{zeros2}
below. Note that the two zeros at $s=\frac{1}{2}$ of $a_0(y,s)$ for 
$y= y^{\ast}$ are omitted from this table. If the zeros
were renumbered to include them as the
first two zeros,  it appears that each  zero of $a_0(y^{\ast},s)$
will be closer to the real axis than the corresponding zero
of $I(y^{\ast}, s)$.

\begin{table}[htbp]
\center
\begin{tabular}{|r|c|c|}
\hline
 & $y=1$ & $y=y^{\ast}$  \\
\hline
1 & 6.974683133 & 2.244794235 \\
2 & 10.40228756 & 3.851296383 \\
3 & 12.42264167 & 5.404657031 \\
4 & 15.08382464 & 6.732441081 \\
5 & 16.40456028 & 7.383718196 \\
6 & 18.68201963 & 8.670185248 \\
7 & 20.34995710 & 10.02271471 \\
8 & 21.60499108 &  10.69728308 \\
9 & 23.85087057 & 11.78575276 \\
10 & 24.83364580 & 12.56610869 \\
11 & 26.40277087 & 13.53142535 \\
12 & 28.11180718 & 14.79167003 \\
13 & 29.54150449 & 15.42550847 \\
14 & 30.39424164 & 16.28902291 \\
15 & 32.41487455 & 16.93621484 \\
\hline
\end{tabular}
\caption{Zeros of constant term $a_0(y,s)$ on critical line}
\label{zeros2}
\end{table}

%
%
%
%
\section{Concluding Remarks}

The  observation of Bump (\cite[p. 6]{BCKV00}) 
mentioned  in
the introduction, concerning the Riemann hypothesis
property holding  for
the Fourier coefficients $a_n(y,s)$ ($n \ne 0$), consisted of a
direct verification, and gave no conceptual explanation why
the Riemann hypothesis holds in these cases. 
Bump observed, more strikingly,
that the truth of the Riemann hypothesis for 
certain Fourier coefficients
of metaplectic Eisenstein series would imply
the Riemann hypothesis for the Riemann zeta function and for 
various Dirichlet $L$-functions. 
The function field analogue 
for metaplectic Eisenstein series was 
unconditionally proved by D. Cardon (\cite{Ca96},\cite{Ca97})
in some cases.
Cardon's proofs also are
direct verifications, making  use of the truth of the 
Riemann hypothesis for curves.
Here our result Theorem~\ref{th13} is again a direct
verification without providing a mechanism. 
A challenging  open question is to find
a conceptual explanation (if there is one)
for the Riemann hypothesis property for 
Fourier coefficients
of such Eisenstein series. 

There also remains the
problem  of conceptually explaining 
the difference in 
behavior of the constant term compared to
the other
Fourier coefficients; why must 
the extra condition $y \ge 1$ be
imposed, and what is the meaning of ``exceptional zeros''? 
Theorem~\ref{th12}
suggests the possibility that  the ``truncated'' 
Fourier coefficients
be considered, for there the Riemann hypothesis does hold
for the ``truncated'' zero-th Fourier coefficient for
$y \ge 1$. This encourages  further study  
of  related integrals.
 
One may ask whether
 the integrated versions of the Fourier coefficients
$a_n(y, s)$ over $\sD_T$ satisfy the Riemann hypothesis. 
The integrals in question are 
\begin{equation}~\label{502}
J_{n}(T, s) := \int \int_{\sD_T} a_n(y,s) \frac{dx dy}{y^2}.
\end{equation}
for integer $n$. In effect we are integrating $E(z,s)$ against the
twisted hyperbolic measure $e^{-2\pi  i n x} \frac{dx dy}{y^2}.$
This question
is related to determining the location of the zeros 
of the special function
\begin{equation}~\label{503}
\tilde{K}_s(y) =: \int_y^{\infty} u^{-\frac{3}{2}} K_s (u) du
\end{equation}
for fixed $y >0$. If for all $y>0$  the zeros of $\tilde{K}_s(y)$ 
lie on the imaginary axis, then the Riemann hypothesis will
hold for all $J_{n}(T, s)$ for all $T > 0$.

Another interesting integral of Eisenstein series 
is the integrated version of the
Fourier coefficients
$a_n(y, s)$ over the entire modular surface $\sF$, i.e.
\begin{equation}~\label{504}
K_{n}(s) := \int \int_{\sF} a_n(y,s) \frac{dx dy}{y^2}. \nonumber 
\end{equation}
The 
observation  of Zagier states  that $K_{0}(s) \equiv 0$, because
$K_{0}(s)$ is an eigenfunction of the non-Euclidean Laplacian.
However this is not the case for the other $K_{n}(s)$,
which are potentially interesting functions.

It might also be of interest to study properties of the
vertical distribution of the zeros of these functions,
particularly as $y$ varies.
Should  some analogue of the GUE property hold for these 
zeros? The results of \cite{La04} on a related problem
suggest that one should expect the zeros of these
functions to  be very smoothly spaced,
with no GUE behavior.

The existence of zeros on the real axis 
for individual Epstein
zeta functions was noted
long ago.
In 1964 Bateman and Grosswald \cite{BG64} gave a criterion 
for individual Epstein zeta functions to have a real zero,
whose main term was  $y^{\ast}= 4\pi e^{-\gamma}$ with a 
very  small error term depending on a parameter $k$. 
In terms of this paper, their  parameter $k=y$, and one of  their results
can be rephrased as saying that each  function $E(z, s)$
with $z=x+iy$ always has a real zero when $y \ge 7.0556$
and never has a real zero when $k\le 7.0554$.
For the constant term
$a_0(y,s)$  obtained the exact cutoff value $y^{\ast}=4\pi e^{-\gamma}$
in Theorem~\ref{th13}, and we note 
that the constant term $a_0(y,s)$ is obtained
by averaging over $x$ of $E(z, s)$, for fixed $y$. 

%
%
%
%
\section{Appendix : Stable and semi-stable lattices}

The rank $n$ zeta function $A_{N, K}(s)$
of a number field $K$ is given in \cite{Weng04}
as an Mellin-type integral over
a moduli space of rank $n$ semi-stable lattices; this 
can be reduced to an integral of an Eisenstein series associated
to the group $PSL(n, \ZZ)$ over a certain subset of the allowable
lattices.

The notion of semistability for lattices
was introduced by Stuhler \cite{Stu76}, in analogy with
a notion of semistability for vector bundles on curves,
as given in Harder and Narasimhan \cite{HN75}.
The  definition was given more generally
for  $\sO_K$-lattices, which
are lattices  having a ring of endomorphisms including
$\sO_K$, the ring of algebraic integers of a number
field $K$. For simplicity we treat here the case $K=\QQ$,
where  $\sO_{\QQ}= \ZZ$, and a $\ZZ$-lattice is
just a lattice. We associate to a rank $r$ lattice
$L$ embedded in $\RR^n$ (with $r \le n$) having  basis
$V_L= \ZZ[\bv_1, ..., \bv_r]$ of row vectors 
$\bv_j = (v_{j,1}, ..., v_{j,n})$, so that $V_L$
is an $n \times r$ matrix, a {\em covolume}
\begin{equation}~\label{A1}
\mbox{Vol(L)} := |\det( V_L V_L^T)|^{\frac{1}{2}}.
\end{equation}
The {\em slope} $s(L)$ 
of $L$ to be
\begin{equation}~\label{A2}
s(L) = \frac{\log( Vol(L))}{dim(L)} = \frac{1}{r} \log( Vol(L)).
\end{equation}

Now suppose that $L$ is an  $n$-dimensional lattice $L$ 
embedded in $\RR^n$.
For
 $1 \le r \le n$ we define the invariants
\begin{equation}~\label{A3}
\kappa_r(L) := \min \{ Vol(L') : ~ L' ~
\mbox{an~rank~r~ sublattice~of}~L\}.
\end{equation}
The minimal slope of rank $r$ sublattices of $L$ is 
given by
\begin{equation}~\label{A4}
s_r(L) := \frac{1}{r}\log \kappa_r(L), ~~ 1 \le r \le n.
\end{equation}
In particular $s_n(L) = \frac{1}{n} \log Vol(L)$.
We also artificially define $s_0(L) =0$.

Next we 
plot the points $(0,0)$ and $\{ (r, s_r(L)): 1 \le r \le n\}$ as
points in the plane, and form 
and their convex hull,
which either forms a polygon,
or in degenerate cases  a line segment. 
The relevance of convexity
is the inequality (Stuhler \cite[Prop.2]{Stu76}) 
\begin{equation}
Vol( L_1 \cap L_2) Vol(L_1 + L_2) \le Vol(L_1)Vol(L_2).
\end{equation}
valid for
two sublattices $L_1, L_2$ of a lattice  $L$, cf.
Grayson~\cite[Theorem 1.12]{Gr84}).
  
Now consider the 
lower envelope of this polygon (the lowest points in it when
intersected with vertical lines), which forms  a finite union of
line segments with different slopes, increasing from left to
right, with vertices occurring
at ranks  $0= r_0 < r_1 < \cdots < r_k=n$.
Grayson \cite[p. 608]{Gr84} terms this envelope  the 
{\em  canonical polygon} of $L$, and  
the  endpoints of these segments $\{ (r_j, s_{r_j}(L) \}$
comprise the {\em canonical vertices}. The canonical
vertices always include
$(0,0)$ and $(n, s(L))$. The key fact about
them, due to Stuhler (cf. Grayson \cite[Theorem 1.18]{Gr84}),
 is that there
is a unique sublattice $L_{r_j}$ in $L$ of rank $r_j$
having slope $ s_{r_j}(L)$, and the set of such lattices
are totally ordered by inclusion. This chain of 
lattices is analogous to 
of the Harder-Narasimhan canonical filtration for vector bundles
on curves. 

The definition of Stuhler \cite{Stu76} is the following.
A $\ZZ$-lattice $L$ is {\em semistable} if the canonical
polygon is a line segment joining $(0,0)$ to  $(n, s(L))$,
and it is {\em stable} if there are no other points 
$(r, s_r(L))$ on this line segment. 
It is {\em unstable} otherwise, i.e. if there are
at least two line segments of different slopes in
the canonical polygon.
For example, the lattice $\ZZ^n$ for $n \ge 2$ is 
semistable but not stable, because it  has $\kappa_j(\ZZ^n) =1$
for $1 \le j \le n$ whence the points determining the 
polygon are $\{ (r,0): ~ 0 \le r \le n\}$.

For the general case of integers $\sO_{K}$ of
an  algebraic number field
$K$ the notion of $\sO_{K}$-semistability is defined considering
only the slopes of those rank $r$ sublattices that 
the action of $\sO_K$ takes into themselves; see 
Stuhler \cite{Stu76}, \cite{Stu77}
or Grayson \cite{Gr84} for details.

Returning to the case of $\ZZ$-lattices, the 
space of $n$-dimensional $\ZZ$-lattices can be identified 
with $GL(n,\ZZ) \backslash GL(n, \RR)$, in
which $GL(n, \RR)$  represents
the set of all bases of the lattice, and the action of $GL(n, \ZZ)$
determines equivalent bases. Each point of this space
can be assigned the property: semistable or unstable. 
The property  of being semi-stable is invariant under 
rotations of $\RR^n$, which is an operation that leaves the
canonical polygon unchanged. The property of semistability
is also preserved  under homotheties
$L \mapsto \lambda L$ for $\lambda > 0$, although this
does change the shape of the canonical polygon to 
$(0.0) \cup \{ (r, s_r(L) + \log \lambda): 1 \le r \le n \}$.
These properties allow 
a notion of semistability
(or unstability) to be unambiguously assigned to points of  
$SL(n, \ZZ) \backslash SL(n, \RR)/ O(n, \RR)$. 

Now we  consider the special case $n=2$. The space
$SL(2, \ZZ) \backslash SL(2, \RR)/ O(2, \RR)$ can
be identified with the upper half-plane
$\HH :=\{ z=x+iy: y>0 \} \equiv SL(2,\RR)/O(2,\RR)$ 
under the action of the modular group $SL(2, \ZZ)$
acting by linear fractional transformations.
The point $z \in \HH$ corresponds to the lattice
$L_{z} = \ZZ[1, z]= \ZZ[(1,0), (x,y)]$, which has $det(L_{z}) = y$.
Now restrict to the standard fundamental domain 
$$
\sF = \{ z=x+iy:~-\frac{1}{2} \le x \le \frac{1}{2},~ |z|^2 \ge 1\}
$$
of $SL(2, \ZZ)\backslash \HH$. The shortest vector in the lattice $L_{z}$ 
then has length $1$, which gives $\kappa_1(L_z)=1$.
The canonical polygon is therefore generated by the
points $(0,0), (1,0)$ and $(2, \frac{1}{2} \log (y^2))$. 
The condition for semistability for $z = x+iy \in \sF$
is that  $y \le 1$,
cf. Grayson \cite[Example 1.25]{Gr84}.
Thus the  points in $\sD_{ss}$ given by
\eqn{113} in \S1 represent the moduli space of
rank $2$ semistable lattices over $\QQ$, quotiented by
the action of homothety and rotations of space.

The {\em rank $r$ vector
bundle $L$-function} $\zeta_{\QQ, r}(s)$ associated
to the rational field $\QQ$ that was introduced by
Lin Weng can be expressed (\cite[p. 8]{we3}) as
$$
\zeta_{\QQ, r}(s)= \frac{r}{2} \pi^{-\frac{rs}{2}} \Gamma(\frac{rs}{2})
\int_{\sM_{\QQ,r}[1]} \left(\sum_{x \in \Lambda \backslash {\bf 0}} 
||\bx||^{-rs}\right) d \mu_1(\Lambda)
$$
The measure $\mu_1(\Lambda)$ is the usual measure on lattices induced
from a (suitably normalized)
Haar measure on $GL(n, \RR)$, whose column vectors represent
a basis of $\Lambda$; equivalence of bases
corresponds to a $GL(n, \ZZ)$
action.  The subset $\sM_{\QQ,r}[1]$ corresponds
to the set of semistable lattices of determinant one, choosing a basis
of positive determinant. The inner sum in the integral above is
an Epstein zeta function of the postive definite
quadratic form in $r$ variables
giving the squared
norm of vectors in the lattice $\Lambda$, which is
$Q(\bx) = \bx^T \sB^T \sB \bx$. In the case $r=2$
this can be identified with the non-holomorphic Eisenstein
series in the paper, and the integral above is simplifiable
to  the integral \eqn{113}.

%
%
%
%

\noindent
Jeffrey C. Lagarias, \\
Dept. of Mathematics, \\
University of Michigan \\
Ann Arbor, MI 48109-1109,\\
USA\\
e-mail address\,:\, lagarias@@umich.edu \\

\noindent
Masatoshi Suzuki \\
Graduate School of Mathematics,\\
Nagoya University,\\
Chikusa-ku, Nagoya 464-8602,\\
Japan\\
e-mail address\,:\,m99009t@@math.nagoya-u.ac.jp

\end{document}